\documentclass{amsart}
\usepackage{amscd}
\usepackage{amssymb}
\include{diagrams}

\newarrow{Corr}{<}{-}{-}{-}{>}

\newcommand{\Z}{{\mathbb Z}}
\newcommand{\Q}{{\mathbb Q}}
\newcommand{\R}{{\mathbb R}}
\newcommand{\fa}{\mathfrak a}
\newcommand{\fA}{\mathfrak A}
\newcommand{\fm}{\mathfrak m}
\newcommand{\fn}{\mathfrak n}
\newcommand{\fp}{\mathfrak p}

\newcommand{\fP}{\mathfrak P}

\newcommand{\cC}{\mathcal C}

\newcommand{\OO}{\mathcal O}
\newcommand{\cT}{\mathcal T}
\newcommand{\cO}{{\mathcal O}}
\newcommand{\cH}{{\mathcal H}}

\newcommand{\Gal}{\operatorname{Gal}\,}
\newcommand{\disc}{\operatorname{disc}\,}
\newcommand{\Cl}{\operatorname{Cl}}
\newcommand{\eps}{\varepsilon}
\newcommand{\la}{\langle}
\newcommand{\ra}{\rangle}
\newcommand{\hra}{\hookrightarrow}
\newcommand{\lra}{\longrightarrow}
\newcommand{\too}{\longmapsto}
\newcommand{\ts}{\textstyle}
\newcommand{\rd}{\operatorname{rd}}

\newcounter{Tc}
\newcounter{Pc}
\newcounter{Lc}

\newtheorem{thm}[Tc]{Theorem}
\newtheorem{prop}[Pc]{Proposition}
\newtheorem{lem}[Lc]{Lemma}
\newtheorem{cor}{Corollary}
\newtheorem*{rem}{Remark}

\title{Ideal class groups of cyclotomic number fields III}
\author{Franz Lemmermeyer}
\address{M\"orikeweg 1 \\ 
         73489 Jagstzell \\
         Germany}
\email{hb3@ix.urz.uni-heidelberg.de}

\subjclass{Primary 11 R 21; Secondary 11 R 29, 11 R 18}
\begin{document}

\begin{abstract}
Using an idea going back to Scholz, we construct unramified
abelian extensions of cyclotomic extensions of number fields.
\end{abstract}

\maketitle

\section{Introduction}

One of the most basic results in algebraic number theory is the
fact that, in every finite extension $K/\Q$ of the rationals, at 
least one prime ramifies, i.e., that $|\disc K| > 1$ except when 
$K = \Q$. This lower bound for the discriminant was conjectured
by Kronecker and first proved by Minkowski, whose geometric
methods gave the more precise bound
$$ |\disc K|  \ge \Big[\Big(\frac \pi4 \Big)^s \frac{n^n}{n!} \Big]^2 $$
for number fields of degree $(K:\Q) = n = r+2s$ with exactly $r$ real 
embeddings. Asymptotically, this shows that the root discriminant
$\rd(K) = |\disc K|^{1/n}$ of a number field satisfies $\rd(K) > e^2$
for large $n$. Blichfeldt \cite{Blich1,Blich} could show that 
$\rd(K) > \pi e$ for large $n$, and that $\rd(K) > 2\pi e^{3/2}$ 
if $K$ is totally real, and these bounds were later improved by
Rogers \cite{Rog1,Rog2}.

One possible approach to answering Furtw\"angler's question whether the 
class field tower always terminates was to show that $\rd(K)$ goes
to infinity with $n$, because fields in the class field tower have
constant root discriminant. In his article \cite{Sch2}, Arnold
Scholz asked how good the asymptotic bounds given above are. He 
denoted by $k_n$ a number field of degree $n$ whose absolute value 
of the discriminant is minimal, and put $E_n = |\disc k_n|^{1/n}$ 
(this is the root discriminant of $k_n$). He observed that the fields 
$\Q(\sqrt[n\,]{2}\,)$ provided the upper bound $E_n < 2n$. 
He conjectured that $\lim E_n/n = 0$ and proved that in fact 
$E_n < (\log n)^2$ for certain $n$ by constructing number fields of 
small discriminant. This was accomplished by studying the ray class 
field $k\{(p)\}$ modulo $(p)$ of $k = \Q(\zeta_{p-1})$.

In this article we will show that $k\{(p)\}$ contains 
$K = \Q(\zeta_{p(p-1)})$, and that $k\{(p)\}/K$ is abelian 
and unramified. In particular we will see that Scholz's 
construction gives a subfield of the Hilbert class field of $K$.

Classically, proofs that certain extensions are unramified
were often done by applying Abhyankar's lemma, which gives
sufficient conditions for killing tame ramification in 
extension fields:

\begin{lem}[Abhyankar's lemma]
Let $L_1/K$ and $L_2/K$ be finite extensions of algebraic number fields,
and let $L = L_1L_2$ denote their compositum. For a prime ideal $\fP$ \
in $L$, let $e_j$ ($j = 1, 2$) denote the ramification indices of
the prime $\fp_j = \fP \cap L_j$ below $\fP$ in $L_j/K$. If 
$e_2 \mid e_1$ and if $\fp_2$ is tamely ramified in $L_2/K$, 
then $\fP$ is unramified in $L/L_1$.
\end{lem}

Abhyankar's Lemma can be used to construct unramified extensions
quite easily: its disadvantage lies in the fact that it cannot handle
wild ramification. We will circumvent this problem by describing the
extensions via ideal groups and using basic class field theory. Our 
main result is the construction of class fields corresponding to 
various factors of class numbers of certain cyclotomic number fields 
that have been found over the last fourty years. Such factors have 
been constructed using 
\begin{itemize}
\item the analytic class number formula (Mets\"ankyl\"a \cite{Met}); 
\item the ambiguous class number formula (Watabe \cite{Wat});
\item Abhyankar's Lemma (Cornell \cite{C,C2,C3},  Madan \cite{Mad},
      Gold \& Madan \cite{GM});
\item Iwasawa theory (Ozaki \cite{Oza});
\item Stickelberger's theorem (Schmidt  \cite{S}).
\end{itemize}

\medskip \noindent
{\bf Remarks.} 
1. The negative solution of the class field tower problem by 
Golod and Shafarevich implies the existence of number fields 
with arbitrarily large degree and bounded root discriminant, 
hence $\liminf_{n \to \infty} E_n/n = 0$ (this follows already
from Scholz's results given above). Scholz's original conjecture 
that $\lim_{n \to \infty} E_n/n = 0$ seems to be still open. 

\medskip\noindent
2. The best upper bounds for $\lim \inf E_n$ used to come from 
examples due to Martinet, whose records have recently been 
improved by Hajir and Maire \cite{HM1,HM2}.

\medskip\noindent
3. Scholz communicated most of the results in \cite{Sch2} to Hasse
in a letter from Aug. 22, 1936 (see \cite{HS}). In this letter he
said he doubted that $E_p = O(\frac{p}{\log p})$ where $p$ runs
through the primes.  

\medskip\noindent
4. There are a lot of open questions regarding the behaviour of
root discriminants. The following problem is particularly appealing 
and might well be accessible with the tools we know today. Let us call 
a group $G$ metabelian of level $m$ if the $m$-th derived group
$G^{(m)} = 1$, but $G^{(m-1)} \ne 1$. Also let $\log_r$ denote the 
$r$-th iterated logarithm, i.e., $\log_0 n = n$, $\log_1 n = \log n$,
and $\log_{r+1} n = \log \log_r n$. Ankeny \cite{Ank} showed that 
there is a constant $c > 0$ such that $\log \rd(K) > c \log_m (K:\Q)$ 
for all normal extensions $K/\Q$ whose Galois group $G$ is metabelian 
of level $m \ge 1$. It seems not to be known whether this is best 
possible in the following sense: for any $m \ge 1$, does there exist 
a sequence of metabelian extensions $K/\Q$ of level $m$ such that 
$\log \rd(K) = O(\log_m (K:\Q))$? The answer to this question 
is clearly positive for $m=1$, where the abelian extensions 
$K = \Q(\zeta_p)$ satisfy 
$\log \rd(K) = \frac{p-2}{p-1} \log p < \log (K:\Q)$. Scholz's 
results imply that the answer is also positive for $m = 2$ since
$\log \rd(K) < 2 \log_2 (K:\Q)$ for the metabelian 
extensions he constructed.

\section{Scholz's construction}

We start by introducing the notation and by explaining the relevant 
facts from class field theory (based on Hasse's exposition \cite{H}). 
Let $k$ be a number field; a {\em modulus} is the symbolic product 
of an integral ideal $\fm$ in $\OO_k$ and various real places; the 
product of all real places of $k$ is denoted by $\infty$. For infinite
places $\fp$, the congruence $\alpha \equiv 1 \bmod \fp$ means that
$\sigma(\alpha) > 0$ for the real embedding $\sigma: k \hra \R$
corresponding to $\fp$.

Below, we will construct ray class fields modulo $\fm = m\OO_k$;
to this end, we need a few definitions:
\begin{itemize}
\item $D_k$ is the group of fractional ideals in $k$;
\item $D_\fm^{\phantom{1}} = \{\fa \in D_k \mid 
        (\fa,\fm) = 1\}$ is the group of ideals $\fa$ coprime to $\fm$;
\item $H_\fm^{\phantom{1}} = \{\fa \in D_\fm \mid \fa = (\alpha) \}$ 
        is its subgroup of principal ideals;
\item $H_\fm^{1} = \{\fa \in D_\fm \mid \fa = (\alpha), 
        \alpha \equiv\ 1 \bmod \fm \}$ is the ray mod $\fm$;
\item $\cH^{\phantom{1}} = \{ \fa \in D_\fm \mid N_{k/\Q}\fa = (a), 
        a > 0, a \equiv 1 \bmod m \}$;
\item $E = E_k = \OO_k^\times$ is the unit group of $\OO_k$;
\item $E_\fm^{1} = \{\eta \in E \mid \eta \equiv\ 1 \bmod \fm \}$.
\end{itemize}

To every ideal group $H_F$ with 
$H_\fm^{1} \subseteq H_F \subseteq D_\fm$ 
class field theory associates a unique abelian extension $F/k$ 
unramified outside $\fm$. The class field corresponding to the 
ideal group $H_\fm^{1}$ is called the {\em ray class field} of $k$ 
modulo $\fm$ and will be denoted by $k\{\fm\}$; the extension 
$k\{\fm\}/k$ is unramified outside $m$ and abelian with Galois 
group $\Gal(k\{\fm\}/k) \simeq D_\fm/H_\fm^{1}$.
The Artin isomorphism shows that $\Gal(F/k) \simeq D_\fm/H_F$, 
and Galois theory then gives $\Gal(k\{\fm\}/F) \simeq H_F/H_\fm^{1}$.

The following result will be of central importance for our construction:

\begin{thm}[Translation Theorem]
If $F$ is the class field of $k$ to the ideal group $H_F$ defined
in $k$, and if $K/k$ is a finite extension, then $FK$ is the class 
field of $K$ to the ideal group 
$$\cT_K(H_F) = \{\fA \in D_\fm(K): N_{K/k}\fA \in H_F\} $$
of all ideals $\fA$ in $K$ coprime to $\fm$ such that $N_{K/k}\fA$ 
is an ideal in $H_F$.
\end{thm}

\begin{figure}[ht!]
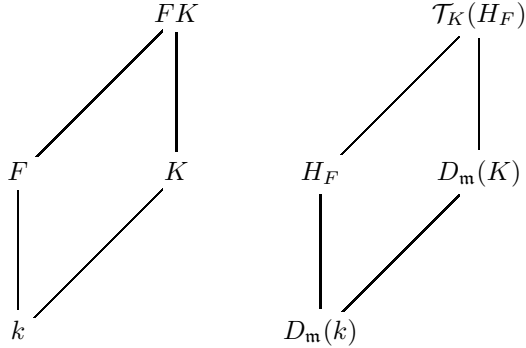
\label{FigVS}
\smallskip
\hbox{\hskip 2cm
\begin{diagram}[p=0.2em]
         &         &   FK \\
         & \ruLine & \dLine \\
      F  &         &  K \\
  \dLine & \ruLine &    \\
      k  &         &    \\
\end{diagram} \hskip 1cm
\begin{diagram}[p=0.2em]
         &         &  \cT_K(H_F) \\
         & \ruLine & \dLine \\
    H_F  &         &  D_\fm(K) \\
  \dLine & \ruLine &  \\
    D_\fm(k)  &         &   \\
\end{diagram}}
\caption{Translation Theorem}
\end{figure}

The Translation Theorem immediately yields

\begin{lem}
The class field of $k$ belonging to the class group $\cH$ is the 
cyclotomic extension $k(\zeta_m)/k$. 
\end{lem}

Since the ideal group attached to the Hilbert class field $k^1/k$
is the group $H_\fm$ of principal ideals defined mod $\fm$ (here
we use the common identification of ideal groups explained in \cite{H}),
and since the intersection of ideal groups corresponds to the 
composita of the class fields, we find the ideal groups attached to the
intermediate fields of $K/k$ displayed in Fig. \ref{Fig1}.

\begin{lem}\label{Lcf}
Let $L$ be the class field of $K$ corresponding to the ideal group $H$
defined mod $\fm$. Then the maximal unramified subextension of $L$ is 
class field for the ideal group $H \cdot H_\fm$, where $H_\fm$ is the
ideal group in $K$ consisting of principal ideals coprime to $\fm$. In 
particular, $L$ is unramified over $K$ if and only if $H_\fm \subseteq H$.
\end{lem}

\begin{proof}
The maximal unramified subextension of $L$ is the intersection 
$L \cap K^1$ of $L$ with the Hilbert class field $K^1$ of $K$;
its associated ideal group is therefore the product of the 
ideal groups $H$ attached to $L$ and $H_\fm$ attached to $K^1$.

Next, $L/K$ is unramified if and only if $L \subseteq K^1$,
that is, if and only if $H \cdot H_\fm = H$; the claim now follows.
\end{proof}

We also need to explain what we mean by the relative class group.
Let $K/k$ be an extension of number fields; the relative norm 
$N_{K/k}$ maps $\Cl(K)$ to $\Cl(k)$, and its kernel $\Cl(K/k)$ 
is called the relative class group of $K/k$. If $K$ is a CM-field 
with maximal real subfield $K^+$, then we put $\Cl^-(K) = \Cl(K/K^+)$.
Finally, if $L/K$ is an extension of CM-fields, then we define
$\Cl^-(L/K)$ to be the intersection of the kernels of the norm maps 
$N_{L/K}:  \Cl(L) \lra \Cl(K)$ and $N_{L/L^+}: \Cl(L) \lra \Cl(L^+)$.

\medskip\noindent{\bf Remark.}
There is a second group measuring the ``new class group'' of $K/k$ 
(which will play no role here), namely the quotient 
$\Cl^*(K/k) = \Cl(K)/\Cl(k)^j$, where $\Cl(k)^j$ is the image of 
the transfer homomorphism $j: \Cl(k) \lra \Cl(K)$. 
\medskip

Now we can state our main theorem:

\begin{thm}\label{SchT} 
Let $k$ be a totally complex number field of degree $(k:\Q) = n$,
and let $p$ be a prime that splits completely in $k/\Q$.
Assume that $m = p^f > 2$, put $\fm = m\OO_k$, and let $w$ denote 
the number of roots of unity contained in $k$. 
Let $M = \Q(\zeta_m)$, $F = kM = k(\zeta_m)$, and $K = k^1(\zeta_m)$, 
where $k^1$ denotes the Hilbert class field of $k$ (see the Hasse
diagram in Fig. \ref{Fig1}).

Then $k\{\fm\}$ is an unramified abelian extension of $F$ 
containing $K$, with relative degree
\begin{equation}\label{CNF}
  (k\{\fm\}:F) \ = \ \frac1{w} h(k)  \phi(m)^{n/2}
                 \cdot \big(E_\fm^{1}:E^{\phi(m)}\big). 
\end{equation}
Moreover, the relative class group $\Cl(F/k)$ contains a subgroup 
$C$ of type
$$ \ts C \simeq \Z/\frac{\phi(m)}{w} \times (\Z/\phi(m))^{\frac n2 -1}.$$ 
If $k$ is a CM field, then so is $F$, and $\Cl^-(F/k)$ contains a 
subgroup of type $C$.
\end{thm}

\begin{proof}
The proof is quite involved, so let me outline the main ideas first.
We show that $F = k(\zeta_m)$ is contained in the ray class 
field $k\{\fm\}$, and that $k\{\fm\}/F$ is unramified. Using
the translation theorem, we can transfer the ideal groups attached
to $k\{\fm\}/F$ to $F$, where they are still defined mod $\fm$.
Since $k\{\fm\}/F$ is unramified, they can be identified with
ideal groups defined mod $(1)$. Finally we observe that these ideal 
groups are killed by various relative norms.

\begin{figure}
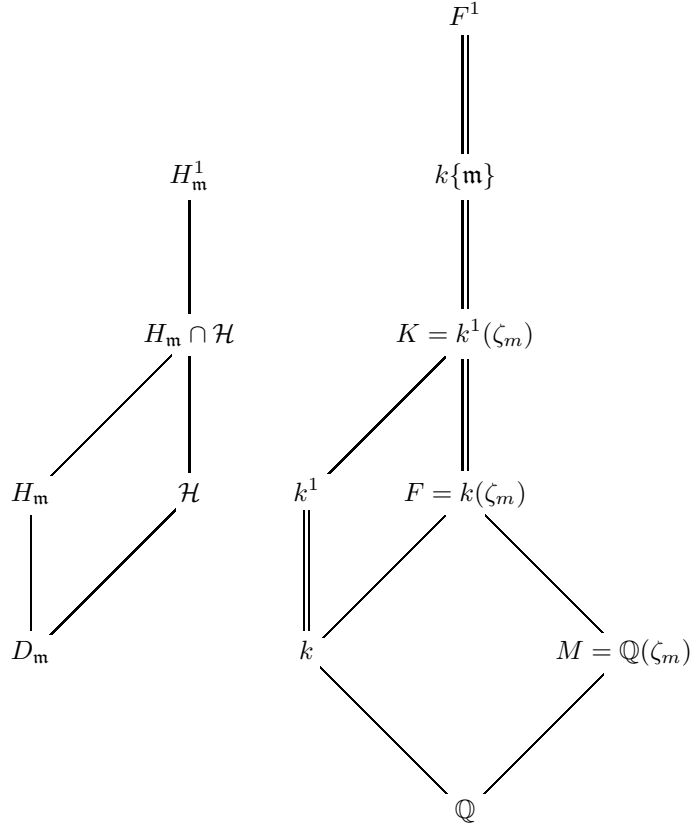

\smallskip
\hbox{\hskip 1.5cm
\begin{diagram}[p=0.2em]
  &  &   \\
  &  &   \\
  &     &   H_\fm^{1} \\
  &     &   \dLine     \\
  &     &   H_\fm\cap \cH \\
  & \ruLine & \dLine \\
  H_\fm &   & \cH \\
  \dLine & \ruLine &  \\
  D_\fm     &       &   \\
    &  &  \\
    &  &  \\
\end{diagram} \hskip 1cm
\begin{diagram}[p=0.2em]
  &     &       F^1       &  & \\
  &     &   \dLine\dLine  &  &  \\
  &     &   k\{\fm\}           &   & \\
  &     &   \dLine\dLine       &   &  \\
  &     &   K = k^1(\zeta_m)   &   & \\
  & \ruLine & \dLine    \dLine &   &  \\
  k^1 &   & F = k(\zeta_m)     &   &  \\
  \dLine\dLine & \ruLine &  &  \rdLine &    \\
  k     &       &  & & M = \Q(\zeta_m)   \\
  & \rdLine &  &  \ruLine &   \\
  &         & \Q &    &
\end{diagram}}
\caption{Ideal Groups and Fields occurring in the Proof of Thm. \ref{SchT}.
A $\parallel$ in the diagram indicates 
that the corresponding extension is unramified.}\label{Fig1}
\end{figure}

\subsection*{1. $F$ is contained in $k\{\fm\}$.}
This follows if we can show that the ideal group $H_\fm^1$ attached
to $k\{\fm\}$ is contained in the ideal group $\cH$ attached to $F$.
The ideal group corresponding to the class field $\Q(\zeta_m)$ over
$\Q$ is $\{(a): a \in \Z, a > 0, a \equiv 1 \bmod m\}$. By the 
translation theorem, this implies that $F$ is the class field for
$\cH_1 = \{\fa \in D_k: N_{k/\Q} \fa = (a), a > 0, a \equiv 1 \bmod m\}$.
Identifying it with the corresponding group $\cH = \cH_1 \cap D_\fm$,
we see that $F$ is the class field for the ideal group 
$\cH : \{\fa \in D_\fm: N_{k/\Q} \fa = (a), a > 0, a \equiv 1 \bmod m\}$.
Now clearly $H_\fm^1 \subseteq \cH$: in fact, if $\fa = (\alpha)$
for $\alpha \equiv 1 \bmod \fm$, then $N_{k/\Q} \fa = (a)$ for 
$a = N_{k/\Q} \alpha$, and clearly $a > 0$ (since $k$ is totally complex)
and $a \equiv 1 \bmod m$ (since $\alpha \equiv 1 \bmod m \cO_k$ implies
$N_{k/\Q} \alpha \equiv 1 \bmod m$).

\subsection*{2. $k\{\fm\}/F$ is unramified.}
By the translation theorem, $k\{\fm\}$ is the class field for the 
ideal group 
$$ H = \{\fA \in D_\fm(F): N_{F/k} \fA 
              = (\alpha), \alpha \equiv 1 \bmod \fm\}. $$
According to Lemma \ref{Lcf}, we need to show that $H$ contains the 
group $H_\fm(F)$ of principal ideals coprime to $\fm$.

To this end, let $A \in F^\times$ be an integer coprime to $\fm = (p^f)$.
We will show that $N_{F/k} A \equiv 1 \bmod \fp^f$ for every prime $\fp$ 
in $k$ above $p$. Since $F/k$ is completely ramified at $\fp$, 
every $A \in F^\times$ coprime to $\fp$ is congruent mod $\fp$ to 
an element in the inertia field, and so we have $A \equiv a \bmod \fp$ 
for some $a \in k$. But $\cO_k/\fp \simeq \Z/p\Z$ since $p$ splits 
completely in $k$, and now it is easy to see that we even have 
$\cO_F/\fp^f \simeq \Z/p^f\Z$. Thus we can choose $a \in \Z$ such 
that $A \equiv a \bmod \fp^f$. 
Taking norms shows that $N_{F/k} A \equiv a^{(F:k)} = a^{\phi(m)} 
\equiv 1 \bmod \fp^f$, which is what we needed to prove.

\medskip\noindent
{\bf Remark.} Let $\xi_m = \zeta_m + \zeta_m^{-1}$ and put 
$F_0 = k(\xi_m)$; then the extension $k\{\fm\}/F_0$ is abelian,
but not unramified. Since $\Gal(F_0/k)$ acts transitively on
half of the prime ideals above $p$ in $F_0$, the compositum of all 
the inertia subgroups has index $\le 4$ in $\Gal(k\{\fm\}/k)$.
This implies that at least a quotient of $C$ by $\Z/2$ already
lives in the class groups $\Cl(F_0/k)$ and $\Cl^-(F_0/k)$, 
respectively.

\subsection*{3. Computation of $(k\{\fm\}:F)$.}
By now we know that the class number of $F$ is divisible
by $(k\{\fm\}:F)$; clearly
$(k\{\fm\}:F) = (k\{\fm\}:k)/\phi(m)$, so it remains to
compute the index $(k\{\fm\}:k)$. The formula for the number of ray
classes gives
$$ (k\{\fm\}:k) = (D_\fm:H_\fm^{1}) 
                     = h(k) \Phi_k(\fm)/(E:E_\fm^{1}).$$
Now $\Phi_k(\fm) = \prod \Phi_k(\fp_j^f)  = \phi(m)^n$; 
moreover, $E^{\phi(m)} \subseteq E_\fm^{1}$, and 
$(E:E^{\phi(m)}) = w \phi(m)^{n/2-1}$. In fact, every 
element of the $\Z$-basis of the unit group $E$ accounts
for a factor of $\phi(m)$, and the unit $\zeta_w$ contributes a 
factor $w$: since $p$ splits completely in $k$, we deduce
that $p \equiv 1 \bmod w$, hence $\zeta_w^{\phi(m)} = 1$. Thus
$$ (E:E_\fm^{1}) = (E:E^{\phi(m)})/(E_\fm^{1}:E^{\phi(m)})
   =  w \phi(m)^{n/2-1}/(E_\fm^{1}:E^{\phi(m)}),$$
and now we find
$$(k\{\fm\}:k) =  \ts{\frac1{w}} h(k)
        \phi(m)^{n/2+1}(E_\fm^{1}:E^{\phi(m)}).$$
Dividing through by $(F:k) = \phi(m)$ we get the
class number factor (\ref{CNF}) in Thm. \ref{SchT}.

\subsection*{4. $C$ is a subgroup of the ray class group.}
Next we are going to realize the group
$$ \ts \Z/\frac{\phi(m)}{w} \times (\Z/\phi(m))^{n/2-1}$$ 
as a subgroup of the ray class group 
$(H_\fm \cap \cH)/H_\fm^{1} \simeq \Gal(k\{\fm\}/F)$.
To this end we first formulate the following

\begin{lem}\label{L3}
Let $F/K$ be an extension of number fields with degree $n$, 
and let $\fm$ be an ideal in $\OO_K$ such that every prime 
dividing $\fm$ splits completely in $F/K$. Then the norm map 
$N = N_{F/K}$ induces an exact sequence
$$ \begin{CD}
  1 @>>> ((\OO_K/\fm)^\times)^{n-1} @>>> (\OO_F/\fm)^\times 
    @>{N}>> (\OO_K/\fm)^\times @>>> 1. \end{CD}$$
\end{lem}

\begin{proof}
Write $\fm = \fp^a \fn$ for some ideal $\fn$ prime to $\fp$,
and $\fp\OO_F = \fP_1 \cdots \fP_r$. Given a residue class 
$\alpha \bmod \fm$ with $\alpha \in \OO_K$, use the Chinese 
remainder theorem to find a $\beta \in \OO_F$ with 
$\beta \equiv \alpha \bmod \fP_1^a$,
$\beta \equiv 1 \bmod \fP_i^a$ for $2 \le i \le r$, and 
$\beta \equiv 1 \bmod \fn$. Then $N_{F/K} \beta \equiv
\alpha \bmod \fm$. This proves that the norm map is onto.

Now for the kernel: since $(\OO_F/\fm\OO_F)^\times \simeq
((\OO_K/\fm\OO_K)^\times)^n$, and since the image of the map
is $(\OO_K/\fm\OO_K)^\times$, it is sufficient to show that
a subgroup isomorphic to $((\OO_K/\fm\OO_K)^\times)^{n-1}$ is in
the kernel. This is done as follows: first observe that we
can write $\fm\OO_F$ as a product $\fm\OO_F = \fm_1 \cdots \fm_n$
of relatively prime conjugate ideals $\fm_i$ such that
$\OO_F/\fm_i \simeq  \OO_K/\fm\OO_K$. Now given any 
$(a_1, \ldots, a_{n-1}) \in (\OO_F/\fm_1)^\times \times \cdots
\times (\OO_F/\fm_{n-1})^\times$, choose algebraic integers 
$\alpha_1, \ldots, \alpha_{n-1} \in \OO_F$ such that 
$\alpha_j \equiv a_j \bmod \fm_j$ and $\alpha_j \equiv 1 \bmod \fm_i$ 
for every $i \ne j$. Now pick $\alpha_n \in \OO_F$ such that 
$\prod_{j=1}^n \alpha_j \equiv 1 \bmod \fm_n$ and
$\alpha_n \equiv 1 \bmod \fm_1\cdots \fm_{n-1}$. Clearly,
the $\alpha_n$ are in the kernel of the norm map, and 
different vectors $(a_1, \ldots, a_{n-1})$ give different
residue classes $\alpha_n \bmod \fm$. This concludes the proof.
\end{proof}

The map sending $\alpha + \fm \in (\OO_k/\fm)^\times$ to the
coset $(\alpha)H_\fm^{1} \in H_\fm/H_\fm^{1}$ induces the
familiar exact sequence in the first row of the diagram
$$ \begin{CD}
 1 @>>> E/E_\fm^{1} @>>> (\OO_k/\fm)^\times @>>> H_\fm/H_\fm^{1} @>>> 1 \\
 @.     @VN_{k/\Q}VV @VN_{k/\Q}VV @VN_{k/\Q}VV   @.   \\
 @.   1 @>>>  (\Z/m)^\times  @>>> (\Z/m)^\times @>>> 1. \end{CD} $$ 
Applying the norm $N_{k/\Q}$ to the groups in the first row we see
that $N_{k/\Q} E = 1$ since $k$ is totally complex; moreover, $N_{k/\Q}$
maps $(\OO_k/\fm)^\times$ to $(\Z/m)^\times$, and Lemma \ref{L3} shows
that the kernel $A_\fm$ of this map is isomorphic to $(\Z/\phi(m))^{n-1}$.

Finally, $N_{k/\Q}$ maps ideals $(\alpha) \in H_\fm$ to ideals
$(a) \in H_{m\infty}$; thus we get a map  
$N_{k/\Q}: H_\fm/H_\fm^{1} \lra H_{m\infty}/H_{m\infty}^{1}$
whose kernel consists of all elements $(\alpha) H_\fm^{1}$
such that $N_{k/\Q} (\alpha) = (a)$ for $a \equiv 1 \bmod m\infty$.
This immediately shows that the kernel is just the ideal group
$\cH \cap H_\fm/H_\fm^{1}$. The group $H_{m\infty}/H_{m\infty}^{1}$
is the ray class group in $\Q$ attached to the cyclotomic extension
$\Q(\zeta_m)/\Q$ and is thus isomorphic to 
$\Gal(\Q(\zeta_m)/\Q) \simeq (\Z/m)^\times$; this can, of course,
also be checked directly by sending an ideal $(a) \in H_{m\infty}$ 
with $a > 0$ to its residue class $a +m\Z$.

The snake lemma then provides us with the exact sequence
\begin{equation}\label{ESF} \begin{CD}
 1 @>>> E/E_\fm^{1} @>>> A_\fm @>>> \cH \cap H_\fm/H_\fm^{1} @>>> 1,
\end{CD} \end{equation}

Clearly $E/E_\fm^{1}$ is a factor group of 
$E/E^{\phi(m)} \simeq \Z/w \times (\Z/\phi(m))^{n/2-1}$; this shows that
$\cH \cap H_\fm/H_\fm^{1}$ contains a subgroup $\cC/H_\fm^{1} 
\simeq C = \Z/\frac{\phi(m)}{w} \times (\Z/\phi(m))^{n/2-1}$ as 
claimed.

\subsection*{5. $C$ is a subgroup of Cl$(F/k)$}

We now use the Translation Theorem to lift the ideal groups 
attached to $K$ and $k\{\fm\}$ from $k$ to $F$, where the
groups will have conductor $1$ since the extensions are unramified.
\begin{enumerate}
\item $\cT(\cH) = D_\fm(F)$.
\item $\cT(H_\fm \cap \cH) = \cT(H_\fm) 
        = \{\fA \in D_\fm(F): N_{F/k} \fA = (\alpha) \}$.
\item $\cT(H_\fm^{1}) = 
        \{\fA \in D_\fm(F): N_{F/k} \fA \in H_{\fm}^{1} \}$.
\end{enumerate}
Since we know that the class field to $\cT(H_\fm)$ is unramified,
it follows that it can be identified with an ideal group $I_F$ defined
modulo $(1)$; recall that $\cT(H_\fm)$ is the group you get by omitting
all ideals not coprime to $\fm$ from $I_F$. Since applying the relative 
norm $N_{F/k}$ to the ideals in $I_F$ gives principal ideals, it is 
immediately clear that the corresponding class group $I_F/H_F$ is 
contained in the relative class group $\Cl(F/k)$.

\subsection*{6. $C$ is a subgroup of the minus class group}
Assume that $k$ is a CM-field, and let $\nu = N_{k/k^+}$ denote 
the relative norm from $k$ to its maximal real subfield $k^+$.
Consider the ideal group $\cH^+ = \{\fa \in D_m: \fa = (\alpha),
 \nu \alpha \equiv 1 \bmod (m)\}$. We claim that there is an exact
sequence
$$ \begin{CD}
  1 @>>> W_k @>>> (\cO/\fm)^\times[\nu] @>{\omega}>> \cH^+/H_\fm^1 @>>> 1,
  \end{CD} $$
where $(\cO/\fm)^\times[\nu]$ is the subgroup of $(\cO/\fm)^\times$
killed by $\nu$. 

The homomorphism 
$\omega: (\cO/\fm)^\times[\nu] \lra \cH^+; \alpha + \fm \too (\alpha)$
is surjective by definition, and its kernel consists of all
residue classes $\alpha + \fm$ with $(\alpha) = (1)$ and 
$\nu \alpha \equiv 1 \bmod m$. The only unit in $k^+$ congruent
to $1 \bmod m$ is $1$, hence the kernel consists of all $\alpha + \fm$
for which $\alpha$ is a unit in $\cO_k^\times$ with relative norm $1$.
The only such units are the roots of unity, and this shows that 
$\ker \omega$ is just the image of $W_k$ in $(\cO/\fm)^\times$.

From Lemma \ref{L3} we get 
$(\cO/\fm)^\times[\nu] \simeq (\cO_{k^+}/(m))^\times
  \simeq ((\Z/m)^{\times})^{n/2}$, and this implies
$\cH^+/H_\fm^1 \simeq C$.

Lifting $\cH^+$ to an ideal group $\cT(\cH^+)$ defined in $F$, 
we get a subgroup of the group of ideals in $F$ that is killed 
by the relative norm $N_{F/F^+}$. Thus $\Cl^-(F)$ contains a 
subgroup of type $C$. Since this subgroup is also contained 
in $\Cl(F/k)$, we deduce that $\Cl^-(F/k)$ contains a subgroup
of type $C$, and this concludes the proof of Theorem \ref{SchT}.
\end{proof}

\medskip
\noindent{\bf Example.} 
Consider the field $K = \Q(\zeta_{155})$. If we take $m = 31$ 
and $k = \Q(\zeta_5)$, then we find that $\Cl^-(K/k)$ contains a 
subgroup of type $\Z/3 \times \Z/30$. If we take the subfield of 
degree $10$ of $\Q(\zeta_{31})$ as $k$ and $m = 5$, then we find 
the subgroup $\Z/2 \times (\Z/4)^4$ of $\Cl^-(F/k)$, where $F = k(\zeta_5)$. 

\section{Applications}

Here we will derive a couple of corollaries from our result.
To this end we need the following observation:

\begin{lem}\label{Lw}
Let $K = \Q(\zeta_n)$ be a cyclotomic number field with 
$n \not\equiv 2 \bmod 4$, and let $p$ be a prime. Then 
the number $w$ of roots of unity contained in the decomposition
field $k$ of $p$ in $K/\Q$ is given by
\begin{equation}\label{Edefw} 
w = \begin{cases} 2  & \text{if} \ p = 2, \\
  	(p-1,n)  & \text{if} \ p > 3 \ \text{and} \ 2 \mid n, \\
	(p-1,2n) & \text{if} \ p > 3 \ \text{and} \ 2 \nmid n. \end{cases}
\end{equation}
\end{lem}

\begin{proof}
Let us dispose of the case $p=2$ first: if $p$ splits completely
in a field containing $\zeta_w$ with $w > 2$, then we must have 
$p \equiv 1 \bmod w$: contradiction. Thus $w = 2$ if $p = 2$.

Now assume that $p$ is odd; the same argument as above shows
that $w \mid p-1$. Moreover, since $k \subseteq K$ we clearly 
have $w \mid n$ or $w \mid 2n$ according as $n$ is even or odd. 
Thus $w \mid (p-1,n)$ if $n$ is even and $w \mid (p-1,2n)$ if 
$n$ is odd. Conversely, if $n$ is even and $w \mid (p-1,n)$, then 
$p$ splits in $\Q(\zeta_w)$ and $\Q(\zeta_w) \subseteq K$, hence 
$\zeta_w \in k$. The argument for odd $n$ is similar.
\end{proof}

This allows us to reformulate Thm. \ref{SchT} somewhat:

\begin{cor}\label{C1}
Let $n \not\equiv 2 \bmod 4$ be an integer, $p$ a prime not dividing $n$, 
and assume that there is no integer $j$ such that $p^j \equiv -1 \bmod n$. 
Let $k$ be the decomposition field of $p$ in $\Q(\zeta_n)$; then the 
relative minus class group $\Cl^-(K/k)$ of $K = k(\zeta_{p^f})$ contains 
a subgroup of type $\Z/\frac uw \times (\Z/u)^{\frac e2 -1}$, where 
$u = \phi(p^f)$, $w$ is defined as in Lemma \ref{Lw}, and 
$e = \phi(n)/o_{n}(p)$, where $o_{n}(p)$ denotes the order of $p$ modulo $n$.
\end{cor}

\begin{proof}
The condition $p^j \not\equiv -1 \bmod n$ is equivalent to the
decomposition field $k$ of $p$ in $\Q(\zeta_n)$ being complex;
we apply Theorem \ref{SchT} to $k(\zeta_n)/k$ and observe
that $(k:\Q) = \phi(n)/o_n(p)$ and that $k$ contains the $w$-th 
roots of unity. 
\end{proof}

As another corollary we get a result due to B. Schmidt 
\cite[Thm. 3.3]{S}; he calls a prime $p$ self-conjugate modulo 
$m = p^an$ if there is an integer $j$ such that 
$p^j \equiv -1 \bmod n$. 

\begin{cor}\label{CSch}
Assume that $p$ is not self-conjugate modulo $m = p^an$. Then 
the minus class group $\Cl^-(K)$ of $K = \Q(\zeta_m)$ has a 
subgroup of type $$ (\Z/w_0\Z) \times (\Z/\phi(p^a)\Z)^{e/2-1},$$
where $e = \phi(n)/o_{n}(p)$, $w_0 = \phi(p^a)/w$, and $w$ is
defined by (\ref{Edefw}).
\end{cor}

\begin{proof}
This follows at once from Corollary \ref{C1}.
\end{proof}

Here is another application of our result:

\begin{prop}\label{P14}
Let $p$ be an odd prime, $L = \Q(\zeta_{4p^2})$, and let $K$ be the
cyclic extension of degree $p$ over $k = \Q(\zeta_4)$ contained in $L$.
Then $p \mid h^-(K)$ if $p \equiv 1 \bmod 4$ and $p \nmid h(K)$
if  $p \equiv 3 \bmod 4$.
\end{prop}

\begin{proof}
It follows from Theorem \ref{SchT} that $p \mid h^-(K)$ if 
$p \equiv 1 \bmod 4$, so assume that $p \equiv 3  \bmod 4$. 
Then only the prime ideal $p$ ramifies in $K/k$, and the 
ambiguous class number formula shows $p \nmid h(K)$. 
\end{proof}

Using the same proof, we can show 

\begin{prop}
Let $p \equiv 1 \bmod 4$ be prime, and let $\fp$ denote a prime
ideal in $k = \Q(i)$ above $p$. Then there exists a 
$\Z_p$-extension $k_\infty/k$ which is unramified outside
$\fp$; moreover, $k_\infty K_\infty$ is an unramified
$\Z_p$-extension of the cyclotomic $\Z_p$-extension $K_\infty$
of $k$.
\end{prop}

It is also possible to prove this by applying the ``blowing up'' 
results from \cite{L2} to the results of Proposition \ref{P14};
one immediately sees that the class group of the $n$-th 
layer of the cyclotomic $\Z_p$-extension of $K$ has a subgroup
isomorphic to $\Z/p^{n+1}\Z$ if $p \equiv 1 \bmod 4$.

\subsection*{Mets\"{a}nkyl\"{a}'s Factorization revisited}
In \cite{L1}, we derived the following factorization of
the class number of certain CM-fields due to 
Mets\"{a}nkyl\"{a} \cite{Met}:

\begin{prop}
Let $L_1 \subseteq \Q(\zeta_m)$ and $L_2 \subseteq \Q(\zeta_n)$ 
be CM-fields, where $m=p^\mu$ and $n=q^\nu$ are prime powers, 
and let $L=L_1L_2$; then
$$ h^-(L) = h^-(L_1)h^-(L_2)T_1T_2,$$
where  $T_1 = h^-(L_1L_2^+)/h^-(L_1)$ and 
$T_2 = h^-(L_2L_1^+)/h^-(L_2)$ are integers.
\end{prop}

Assume that $L = \Q(\zeta_{pq})$ with $p \equiv 1 \bmod q$; by
what we have shown, $\Cl^-(L)$ contains a subgroup of type
$C = \Z/\frac{p-1}2 \times (\Z/p-1)^{(q-3)/2}$. Is this a factor 
of $T_2$? The Remark at the end of Section 2 in the proof of 
Thm. \ref{SchT} shows that $\Cl^-(L_2L_1^+/L_2)$ contains at least 
a quotient of $C$ by $\Z/2\Z$.

\subsection*{Remarks}
Let me add a few remarks concerning \cite{L1}. 
\cite[Prop. 2]{L1}, which I credited to Louboutin, 
actually already appears as Proposition 3 in \cite{YH1}.
Also, Schoof \cite{SCM} introduced a unit index 
$[\mu_K^{\phantom{-}}:\mu_K^{-}]$ which coincides with 
$2/Q^*$ in Hasse's work.

\section*{Acknowledgement}
I thank the referee for reading the manuscript very carefully.

\end{document}